# TIME CHANGES OF SYMMETRIC DIFFUSIONS AND FELLER MEASURES

By Masatoshi Fukushima[1], Ping He and Jiangang Ying[2]

*Kansai University, SHUFE and Fudan University*

We extend the classical Douglas integral, which expresses the Dirichlet integral of a harmonic function on the unit disk in terms of its value on boundary, to the case of conservative symmetric diffusion in terms of Feller measure, by using the approach of time change of Markov processes.

**1. Introduction.** In the present work, we are concerned with a formula that goes back to Douglas [6]:

$$(1.1) \quad \tfrac{1}{2}\int_D |\nabla Hf(x)|^2\, dx = \tfrac{1}{2}\int_{\partial D\times \partial D\setminus d}(f(\xi)-f(\eta))^2\, U(\xi,\eta)\,d\xi\,d\eta,$$

where $Hf$ denotes the harmonic function on the planar unit disk $D$ with boundary value $f$ and $U(\xi,\eta)=\frac{1}{4\pi(1-\cos(\xi-\eta))}$.

In 1962, Doob [4] extended formula (1.1) to the case where $D$ is a general Green space and $\partial D$ is its Martin boundary by adopting the Naim kernel as $U$. Fukushima [8] identified the Naim kernel with the Feller kernel soon after and then utilized the resulting formula (1.1) as a basis to describe all possible symmetric Markovian extensions of the absorbing Brownian motion on a bounded Euclidean domain [9]. The Feller kernel was introduced by Feller [7] for the minimal Markov process on a countable state space for the purpose of describing all possible boundary conditions on some ideal boundaries. A common feature of the mentioned approaches is that we are only given a minimal process on $D$ a priori and we try to capture its Markovian extensions including the construction of intrinsic boundaries.

Received March 2003; revised October 2003.
[1]Supported by Grand-in-Aid for Scientific Research 1164042 from MEXT.
[2]Supported in part by funding from NSFC Grant 10271109.

Since the sixties, investigations of Markov processes and associated Dirichlet forms have been developed considerably and we can now take the following different but much more stochastic view of the formula (1.1). What is given in advance is the reflecting Brownian motion $X$ on $\overline{D}$ and we consider its time-changed process $Y$ on $\partial D$ with respect to a local time on $\partial D$. The left-hand side of (1.1) is the Dirichlet form for $Y$ (the trace of the Dirichlet form for $X$), while the right-hand side is its specific Beurling–Deny representation. Equation (1.1) tells us that $Y$ is of pure jump and that its jumping mechanism, namely, the Lévy system, is governed by the Feller kernel $U$ which can be easily and intrinsically defined depending only on the absorbing Brownian motion $X_D$ on $D$.

This viewpoint allows us to extend the formula (1.1) with great generality. Indeed, we consider in this paper a general symmetric diffusion process $X$ with a general state space $E$ and its time-changed process $Y$ on an arbitrary closed subset $F$ of $E$. We show in Section 5 that the jumping measure and the killing measure for $Y$ can be identified with the Feller measure $U$ and the supplementary Feller measure $V$, respectively, introduced in Section 2, depending only on the absorbed process $X_G$ on $G = E \setminus F$.

The organization of the present paper is as follows. In Section 2, we consider a conservative Borel right process $X$ which is symmetric with respect to a $\sigma$-finite measure $m$ on a general state space $E$. Let $F$ be a closed subset of $E$ and let $X_G$ be the process on $G = E \setminus F$ obtained from $X$ by killing on its hitting time of $F$. We introduce a notion of the energy functional $L_G$ for each pair of $X_G$-almost excessive functions on $G$, a variant of the concept described by Meyer (see [13]). By means of $L_G$, we can readily define the Feller measure $U(d\xi, d\eta)$ (a bi-measure on $F \times F$), the supplementary Feller measure $V(d\xi)$ (a measure on $F$ informally called an escape measure) and the Feller kernel $U(\xi, \eta)$ when the Poisson kernel exists. In Example 2.1, we exhibit explicit expressions of these quantities for the case that $X$ is the $d$-dimensional Brownian motion ($d \geq 3$) and $F$ is the $(d-1)$-dimensional compact smooth hypersurface.

From Section 3 on, we assume that $X$ is a diffusion, namely, its sample paths are continuous. In Section 3, we focus our attention on excursions of the sample paths of $X$ away from the closed set $F$, and we identify the Feller measure and supplementary Feller measure with the expectations of certain homogeneous random measures generated by the endpoints of excursions. We make use of a description of the joint distribution of endpoints of excursions previously studied by Hsu [14] for reflecting Brownian motion on a smooth domain.

From Section 4 on, we further assume that $X$ is associated with a regular Dirichlet space $(E, m, \mathcal{F}, \mathcal{E})$ (without loss of generality owing to the transfer method). In Section 4, we first prove that $F$ always admits an admissible measure $\mu$ in the sense that $\mu$ charges no set of zero capacity and possesses full quasisupport $F$. We then show, by applying a general reduction theorem formulated in the final section (Section 8), that the time-changed process $Y$ of $X$ with respect to the positive continuous additive functional with Revuz measure $\mu$ can be restricted outside some $X$- and $Y$-polar set to be a Hunt process. This reduction enables us to use a general theorem [11] directly to express the jumping measure and the killing measure in the Beuring–Deny representation of the Dirichlet form for $Y$ by means of the Lévy system of $Y$.

By making use of the results in Sections 3 and 4, we prove in Section 5 the stated main assertion (Theorem 5.1) that the jumping measure and the killing measure for the time-changed process $Y$ are identical to the Feller measure and the supplementary Feller measure, respectively. Theorem 5.1 tells us that the trace Dirichlet form $\mathcal{E}(Hf, Hf)$ always dominates the generalized Douglas integral with the Feller measure. By assuming that $m(G)$ is finite, we prove in Section 6 that they are identical under the condition that the energy measures $\mu_{\langle u \rangle}$ do not charge the set $F$ for any $u \in \mathcal{F}$. This condition is satisfied when the energy measures are absolutely continuous with respect to $m$ [the densities are so-called square field operators $\Gamma(u)$] and $m(F) = 0$. We also characterize this condition in terms of the notion of the reflected Dirichlet space of the part of $\mathcal{E}$ on the set $G$ formulated by Silverstein [19, 20] and Chen [3]. In the course of the proof, we make full use of several results in [11] to recover and extend the method in [4] and [8] for computing the Dirichlet norm of the classical harmonic function.

In Section 7, we apply the obtained results to the reflecting Brownian motion on the closure of a bounded Lipschitz domain $D \subset \mathbb{R}^d$ associated with the Dirichlet space $H^1(D)$. In this case, the relative boundary $\partial D$ is known to be identical with the Martin boundary of $D$, so that Doob's representation of (1.1) is recovered by the present approach.

In Section 8, we formulate a general theorem of reduction of a right process to a Hunt process properly associated with a regular Dirichlet form.

**2. Feller measure $U$, supplementary Feller measure $V$ and Feller kernel.** Let $E$ be a Lusin topological space and let $m$ be a $\sigma$-finite positive Borel measure on $E$. Let $X = (X_t, P^x)$ be a conservative Borel right Markov process on $E$ which is $m$-*symmetric* in the sense that its transition function $p_t$ satisfies

$$\int_E p_t f(x) \, g(x) m(dx) = \int_E f(x) \, p_t g(x) m(dx) \qquad \forall f, g \in \mathcal{B}^+.$$

Fix a closed set $F$ and put $G := F^c$. Denote by $T$ the hitting time of $F$. Let

$$p_t^0(x, A) := P^x(X_t \in A, t < T), \qquad x \in G, \ A \subset G,$$

be the transition function of $X_G$, the absorbed process of $X$ on $G$, which is obtained by killing $X$ on leaving $G$. Then $X_G$ is symmetric with respect to the measure $m_G = \mathbb{1}_G \cdot m$ [11]. The resolvent of $X_G$ is denoted by $R_\alpha^0$.

A measurable function $u$ on $G$ is said to be $\alpha$-excessive for $X_G$ if for every $x \in G$,

$$u(x) \geq 0, \qquad e^{-\alpha t} p_t^0 u(x) \uparrow u(x), \ t \downarrow 0.$$

If the above properties hold for $m_G$-a.e. $x \in G$, then $u$ is said to be $\alpha$-almost-excessive. A 0-excessive (resp. 0-almost-excessive) function is simply called excessive (resp. almost excessive). Let us denote by $\mathcal{S}_G$ the totality of $X_G$-almost-excessive functions on $G$ finite $m_G$-a.e. and let $\langle u, v \rangle_{m_G}$ denote the integral of $uv$ with the measure $m_G$.

LEMMA 2.1. *For any $u, v \in \mathcal{S}_G$,*

$$(2.1) \qquad \frac{1}{t} \langle u - p_t^0 u, v \rangle_{m_G}$$

*is nondecreasing as $t \downarrow 0$. If moreover $v$ is $p_t^0$-invariant in the sense that $p_t^0 v = v$, $t > 0$, then (2.1) is independent of $t > 0$.*

PROOF. We set

$$e(t) = \langle u - p_t^0 u, v \rangle_{m_G}.$$

Then, for $t, s \geq 0$,

$$e(t+s) = e(t) + \langle p_t^0 u - p_{t+s}^0 u, v \rangle_{m_G} = e(t) + \langle u - p_s^0 u, p_t^0 v \rangle_{m_G} \leq e(t) + e(s),$$

the last inequality being replaced by equality if $v$ is $p_t^0$-invariant. □

Let us define the *energy functional* of $u, v \in \mathcal{S}_G$ by

$$(2.2) \qquad L_G(u, v) = \lim_{t \downarrow 0} \frac{1}{t} \langle u - p_t^0 u, v \rangle_{m_G}.$$

We note that $L_G(u, v)$ is nothing but the value of the energy functional of the excessive measure $u \cdot m_G$ and the excessive function $v$ for $X_G$ in the sense of Dellacherie-Meyer and Getoor when $X_G$ is transient and $u \cdot m_G$ is purely excessive ([13], Proposition 3.6). We also have the formula

$$(2.3) \qquad L_G(u, v) = \lim_{\alpha \to \infty} \alpha \langle u - \alpha R_\alpha^0 u, v \rangle_{m_G}$$

as an increasing limit because, by the Fubini theorem,

$$(2.4) \qquad \alpha \langle u - \alpha R_\alpha^0 u, v \rangle_{m_G} = \int_0^\infty e^{-t} (t/\alpha)^{-1} \langle u - p_{t/\alpha}^0 u, v \rangle_{m_G} t \, dt.$$

For $\alpha \geq 0$, let $H^\alpha$ be the $\alpha$-order hitting measure for $F$, that is,

$$H^\alpha(x, B) := E^x(e^{-\alpha T}\mathbb{1}_B(X_T); T < \infty), \qquad x \in G, B \in \mathcal{B}(E).$$

Let $H^0$ be denoted by $H$ and let $H^\alpha(x, \cdot)$ be carried by $F$, since $F$ is closed. It is easy to see that, for any $f \in \mathcal{B}(F)^+$, $H^\alpha f$ is $\alpha$-excessive for $X_G$.

We also consider the function on $G$ defined by

(2.5) $$q(x) = P^x(T = \infty)(= 1 - H1(x)), \qquad x \in G.$$

Then $q$ is not only excessive for $X_G$, but also $p_t^0$-invariant.

We now let, for $f, g \in b\mathcal{B}(F)^+$,

(2.6) $$U(f \otimes g) = L_G(Hf, Hg), \qquad V(f) = L_G(Hf, q).$$

We call $U$ the *Feller measure* for $F$ with respect to $m$ because it is a bimeasure in the sense that $U(I_B \otimes I_C)$ is a (possibly infinite) measure in $B \in \mathcal{B}(F)$ [resp. $C \in \mathcal{B}(F)$] for each fixed $C$ (resp. $B$). Measure $V$ is a (possibly infinite) measure on $F$ and is called the *supplementary Feller measure* or, more informally, the *escape measure* for $F$. We see in Section 4 that $U$ is a $\sigma$-finite measure on $F \times F$ off the diagonal and that $V$ is a $\sigma$-finite measure on $F$.

For $\alpha > 0$, we also define the $\alpha$-order Feller measure $U_\alpha$ for $F$ by

(2.7) $$U_\alpha(f \otimes g) = \alpha \langle H^\alpha f, Hg\rangle_{m_G}, \qquad f, g \in b\mathcal{B}(F)^+.$$

LEMMA 2.2. *The formulae for $f, g \in b\mathcal{B}(F)^+$ are*

(2.8) $$U(f \otimes g) = \lim_{t \to 0} \frac{1}{t} E^{Hg \cdot m_G}(T \leq t, f(X_T)),$$

(2.9) $$U(f \otimes g) = \lim_{\alpha \to \infty} U_\alpha(f \otimes g).$$

PROOF. The first formula follows from

$$P^x(T \leq t, f(X_T)) = Hf(x) - p_t^0 Hf(x), \qquad x \in G.$$

The second formula is a consequence of (2.3) and $H^\alpha f = Hf - \alpha R_\alpha^0 Hf$. □

The notion $U$ goes back to Feller [7], who introduced a version of $U$ by (2.9) and utilized it to describe possible boundary conditions for a minimal Markov process on a countable state space.

The supplementary Feller measure $V$ has more specific properties:

LEMMA 2.3. (i) *For any $t > 0$, $\alpha > 0$,*

(2.10) $$V(f) = \frac{1}{t} E^{q \cdot m_G}(T \leq t, f(X_T)), \qquad f \in b\mathcal{B}(F)^+,$$

(2.11) $$V(f) = \alpha \langle H^\alpha f, q\rangle_{m_G}, \qquad f \in b\mathcal{B}(F)^+.$$

(ii) *If $m(G) < \infty$, then $V = 0$.*

(iii) *If $m(G) < \infty$ and $P^x(T < \infty) > 0$ for $m$-a.e. $x \in G$, then $P^x(T < \infty) = 1$ for q.e. $x \in G$.*

PROOF. Part (i) follows from $p_t^0$ invariance of $q$, Lemma 2.1 and (2.4). If $m(G)$ is finite, then the right-hand side of (2.10) tends to zero as $t \to \infty$ and we get (ii). Part (iii) follows from (i) and (ii). □

When the hitting measure $H(x, \cdot)$ has a suitable density with respect to a certain measure $\mu$ on $F$, then the Feller measure $U$ has also a density with respect to $\mu \times \mu$. In the rest of this section, we assume that there exist a $\sigma$-finite measure $\mu$ on $F$ and a finite-valued function $K(x,\xi)$, $x \in G$, $\xi \in F$, strictly positive $(m_G \times \mu)$-a.e. such that

$$(2.12) \quad H(x,B) = \int_B K(x,\xi) \mu(d\xi) \qquad \forall B \in \mathcal{B}(F) \text{ for } m_G\text{-a.e. } x \in G,$$

and $K(\cdot, \xi)$ is $X_G$-almost excessive for every $\xi \in F$. The function $K^\xi(x) = K(x,\xi)$ is called a *Poisson kernel* with respect to $\mu$.

We put

$$(2.13) \qquad U(\xi,\eta) = L_G(K^\xi, K^\eta), \qquad \xi, \eta \in F,$$

which is called a *Feller kernel* with respect to $\mu$.

In fact, if we define the $\alpha$-order Poisson kernel by

$$(2.14) \qquad K_\alpha(x,\xi) = K(x,\xi) - \alpha R_\alpha^0 K^\xi(x), \qquad x \in G, \ \xi \in F,$$

and the $\alpha$-order Feller kernel by

$$(2.15) \qquad U_\alpha(\xi,\eta) = \alpha \langle K_\alpha^\xi, K^\eta \rangle_{m_G}, \qquad \xi, \eta \in F,$$

then, by (2.3),

$$(2.16) \qquad U(\xi,\eta) = \lim_{\alpha \to \infty} U_\alpha(\xi,\eta), \qquad \xi, \eta \in F$$

and we get from (2.9) that

$$U(d\xi, d\eta) = U(\xi,\eta) \mu(d\xi) \mu(d\eta).$$

EXAMPLE 2.1 (Brownian motion and a compact hypersurface). Let $X$ be the standard Brownian motion on $\mathbb{R}^d$ with $d \geq 3$. Let $S$ be a $C^3$ compact hypersurface so that $G = \mathbb{R}^d \setminus S$ is the union of the interior domain $D_i$ and the exterior domain $D_e$. The absorbed Brownian motion $X_G$ has the transition density

$$(2.17) \quad \begin{aligned} p_t^0(x,y) &= n(t, x-y) - E^x[n(t-T, X_T - y); T < t], \\ n(t,x) &= \frac{1}{(2\pi t)^{d/2}} \exp\left(-\frac{|x|^2}{2t}\right), \end{aligned}$$

where $T$ is the hitting time of $S$ by $X$. Density $p_t^0(x,y)$, $x,y \in D_i$ (resp. $x,y \in D_e$) is the fundamental solution of the heat equation

$$(2.18) \quad \frac{\partial u(t,x)}{\partial t} = \frac{1}{2}\Delta_x u(t,x), \qquad t > 0, \ x \in D_i \text{ (resp. } x \in D_e\text{)},$$

with the Dirichlet boundary condition

$$u(t,x) = 0, \qquad x \in S.$$

Denote by $\sigma$ the surface measure on $S$. Then we can get the expressions

$$P^x(T \in ds, X_T \in d\xi) = g(s,x,\xi)\, ds\, \sigma(d\xi)$$

with

$$(2.19) \qquad g(s,x,\xi) = \frac{1}{2}\frac{\partial p_s^0(x,\xi)}{\partial n_\xi^i}, \qquad x \in D_i, \ \xi \in S,$$

$$(2.20) \qquad g(s,x,\xi) = \frac{1}{2}\frac{\partial p_s^0(x,\xi)}{\partial n_\xi^e}, \qquad x \in D_e, \ \xi \in S,$$

where $n_\xi^i$ and $n_\xi^e$ denote the inward normal and outward normal at $\xi \in S$, respectively. A proof of (2.19) was given in [1], page 262. We give a similar proof of (2.20) for completeness.

We extend a smooth function $h$ on $S$ to $D_e$ by

$$h(y) = E^y(h(X_T); T < \infty), \qquad y \in D_e.$$

It can be seen that $h$ is a harmonic function on $D_e$ vanishing at $\infty$ and hence the first derivative of $h$ are bounded on $D_e$ [see the paragraph below (2.25)]. On the other hand, we can see from (2.17) that, for each $T > 0$ and $a > 0$, there are positive constants $C_1, C_2$ such that

$$(2.21) \qquad \left|\frac{\partial p_t^0(x,y)}{\partial y_k}\right| \leq C_1 \exp(-C_2|x-y|^2),$$

$$0 < t < T, \ 1 \leq k \leq d, \ |x-y| > a.$$

For large $R > 0$, we put $D_e^R = \{x \in D_e : |x| < R\}$ and denote its outer boundary by $\Sigma_R$. For a fixed $x \in D_e^R$, we have by Green's formula

$$\frac{1}{2}\int_S \frac{\partial p_s^0(x,y)}{\partial n_y^e} h(y)\, d\sigma(y)$$

$$= \int_{D_e^R}\left(-\frac{1}{2}\Delta_y p_s^0(x,y)\right) h(y)\, dy$$

$$- \frac{1}{2}\int_{\Sigma_R} \frac{\partial p_s^0(x,y)}{\partial n_y} h(y)\, d\sigma(y) + \frac{1}{2}\int_{\Sigma_R} p_s^0(x,y)\frac{\partial h(y)}{\partial n_y}\, d\sigma(y).$$

By the above observations, the last two integrals vanish as $R \to \infty$. Substituting (2.18) into the resulting equality and integrating in $s$, we arrive at

$$\int_0^t ds \int \frac{1}{2} \frac{\partial p_s^0(x,y)}{\partial n_y^e} h(y) \, d\sigma(y) = h(x) - p_s^0 h(x) = E^x(h(X_t); T \leq t),$$

proving (2.20).

Accordingly, the Poisson kernel and the $\alpha$-order Poisson kernel with respect to $\sigma$ admit the expressions

(2.22)
$$K(x, \xi) = \int_0^\infty g(s, x, \xi) \, ds,$$

$$K_\alpha(x, \xi) = \int_0^\infty e^{-\alpha s} g(s, x, \xi) \, ds, \qquad x \in G, \ \xi \in S.$$

The $\alpha$-order Feller kernel $U_\alpha(\xi, \eta)$ is the sum of $U_\alpha^i(\xi, \eta)$ and $U_\alpha^e(\xi, \eta)$, where

$$U_\alpha^i(\xi, \eta) = \alpha \int_{D^i} K_\alpha^\xi(x) K^\eta(x) \, dx,$$

$$U_\alpha^e(\xi, \eta) = \alpha \int_{D^e} K_\alpha^\xi(x) K^\eta(x) \, dx, \qquad \xi, \eta \in S.$$

From (2.19), (2.20) and (2.22), we can get, for $\xi, \eta \in \partial D$, $\xi \neq \eta$,

$$U_\alpha^i(\xi, \eta) = \frac{1}{4} \int_0^\infty (1 - e^{-\alpha t}) \frac{\partial^2 p_t^0(\xi, \eta)}{\partial n_\xi^i \partial n_\eta^i} \, dt,$$

$$U_\alpha^e(\xi, \eta) = \frac{1}{4} \int_0^\infty (1 - e^{-\alpha t}) \frac{\partial^2 p_t^0(\xi, \eta)}{\partial n_\xi^e \partial n_\eta^e} \, dt.$$

By letting $\alpha \to \infty$, we are led to the following expressions of the Feller kernel:

(2.23) $\quad U(\xi, \eta) = \dfrac{1}{4} \int_0^\infty \dfrac{\partial^2 p_t^0(\xi, \eta)}{\partial n_\xi^i \partial n_\eta^i} \, dt + \dfrac{1}{4} \int_0^\infty \dfrac{\partial^2 p_t^0(\xi, \eta)}{\partial n_\xi^e \partial n_\eta^e} \, dt, \qquad \xi \neq \eta,$

(2.24) $\quad U(\xi, \eta) = \dfrac{1}{2} \dfrac{\partial K(\xi, \eta)}{\partial n_\xi^i} + \dfrac{1}{2} \dfrac{\partial K(\xi, \eta)}{\partial n_\xi^e}, \qquad \xi \neq \eta.$

We consider the special case that $S = \Sigma_R$, the sphere of radius $R$ centered at the origin. The Poisson kernel with respect to the surface measure $\sigma$ is then expressed as

(2.25) $\quad K(x, \eta) = \begin{cases} \dfrac{1}{\Omega_d R} \cdot \dfrac{R^2 - |x|^2}{|x - \eta|^d}, & |x| < R, \ \eta \in \Sigma_R, \\ \dfrac{1}{\Omega_d R} \cdot \dfrac{|x|^2 - R^2}{|x - \eta|^d}, & |x| > R, \ \eta \in \Sigma_R, \end{cases}$

where $\Omega_d$ denotes the area of the unit sphere in $\mathbb{R}^d$. Note that, for $D_e = \{|x| > R\}$ and a continuous function $f$ on $\Sigma_R$,

$$(Hf)(x) = \int_{\Sigma_R} K(x, \eta) f(\eta) \, d\sigma(\eta), \qquad x \in D_e,$$

is the unique harmonic function on $D_e$ taking value $f$ on $\Sigma_R$ and vanishing at $\infty$.

By (2.24), we obtain an explicit expression of the Feller kernel:

$$(2.26) \qquad U(\xi, \eta) = \frac{2}{\Omega_d} |\xi - \eta|^{-d}, \qquad \xi, \eta \in \Sigma_R, \ \xi \neq \eta.$$

We can also obtain an explicit expression of the supplementary Feller measure $V$. By virtue of the above observation, $H1(x) = \frac{R^{d-2}}{|x|^{d-2}}$, $x \in D_e$, and consequently we get from (2.10) and (2.20) that

$$V(d\xi) = v(\xi) \sigma(d\xi)$$

with

$$v(\xi) = \frac{1}{2t} \int_0^t ds \int_{\{|x|>R\}} \left(1 - \frac{R^{d-2}}{|x|^{d-2}}\right) \frac{\partial p_s^0(x, \xi)}{\partial n_\xi^e} \, dx, \qquad \xi \in \Sigma_R.$$

The integral on the right-hand side converges in view of (2.21). This expression shows that $v(\xi)$ is actually a positive constant, say, $v_0$ independent of $\xi$ so that

$$(2.27) \qquad V(d\xi) = v_0 \, \sigma(d\xi), \qquad d\xi \in \mathcal{B}(\Sigma_R).$$

Some computations similar to those above were carried out in [17] for a certain Markov process and also in [8] and [14] for diffusions on an interior Euclidean domain.

**3. Endpoints of excursions and $U$ and $V$.** In the sequel, we further assume that $X$ is a diffusion, namely, all of its sample paths are continuous on $[0, \infty)$. For any $\omega \in \Omega$, we define

$$J(\omega) = \{t \in [0, \infty) : X_t(\omega) \in G\},$$

which is open and consists of all of excursions away from $F$ of the sample path of $\omega$.

We set, for $t \geq 0$,

$$R(t) = \inf(t, \infty) \cap J^c = \inf\{s > t : X_s \in F\}, \qquad \inf \varnothing = \infty,$$

and, for $t > 0$,

$$L(t) = \sup[0, t) \cap J^c = \sup\{0 < s < t : X_s \in F\}, \qquad \sup \varnothing = 0.$$

Clearly $R(t) = T \circ \theta_t + t$ and for any $s, t \geq 0$, $R(t) \circ \theta_s + s = R(t+s)$. By continuity of paths, $X_{R(t)} \in F$ if $R(t) < \infty$ and $X_{L(t)} \in F$ on $T < t$. The process $X_{L(t)}$ stays on $F$ until $X$ hits $F$ again and is adapted, but $X_{R(t)}$ is not adapted in general.

For $t > 0$, we introduce the time reversal operator at $t$ by

$$r_t \omega(s) = \omega(t-s),$$

so that $X_s \circ r_t = X_{t-s}$, $s \in [0, t]$.

Since $X$ is $m$-symmetric and conservative, we have

(3.1) $$E^m(Y \circ r_t) = E^m(Y)$$

for any $\mathcal{F}_t$-measurable random variable $Y$ (cf. [11], Lemma 4.1.2).

We can see that

$$L(t) \circ r_t = t - T, \qquad X_{L(t)} \circ r_t = X_{t-L(t)\circ r_t} = X_T$$

on $T < t$.

LEMMA 3.1. *For $t > 0$, let $I_1 \subset [0, t]$, $I_2 \subset [t, \infty)$ be nonempty intervals and let $A, B \in \mathcal{B}(F)$. Then*

(3.2)
$$P^m(L(t) \in I_1, X_{L(t)} \in A, X_t \in G, R(t) \in I_2, X_{R(t)} \in B)$$
$$= \int_G P^x(T \in t - I_1, X_T \in A) P^x(T \in I_2 - t, X_T \in B) m(dx).$$

*In particular,*

(3.3)
$$P^m(L(t) \in I_1, X_{L(t)} \in A, X_t \in G, X_{R(t)} \in B, R(t) < \infty)$$
$$= \int_G P^x(T \in t - I_1, X_T \in A) H(x, B) m(dx).$$

*Furthermore,*

(3.4)
$$P^m(L(t) \in I_1, X_{L(t)} \in A, X_t \in G, R(t) = \infty)$$
$$= \int_G P^x(T \in t - I_1, X_T \in A) q(x) m(dx).$$

PROOF. Clearly $\{L(t) \in I_1, X_{L(t)} \in A\} \in \mathcal{F}_t$. Since $R(t) = T \circ \theta_t + t$, by the Markov property and (3.1), we have

$$P^m(L(t) \in I_1, X_{L(t)} \in A, X_t \in G, R(t) \in I_2, X_{R(t)} \in B)$$
$$= E^m(L(t) \in I_1, X_{L(t)} \in A, P^{X_t}(X_0 \in G, X_T \in B, T \in I_2 - t))$$
$$= E^m[(\mathbb{1}_{\{L(t) \in I_1, X_{L(t)} \in A\}} \phi(X_t)) \circ r_t]$$
$$= E^m[\phi(X_0); X_T \in A, T \in t - I_1],$$

where $\phi(x) = \mathbb{1}_G(x)P^x(X_T \in B, T \in I_2 - t)$. This completes the proof of (3.2) and (3.3). Equation (3.4) follows from (3.3). □

Denote by $I$ the set of all left endpoints of open (excursion) intervals in $J$. We note for $s > 0$ that $s \in I$ if and only if $R(s-) < R(s)$ and that, in this case, $R(s-) = s$. It is convenient to add an extra point $\Delta$ to $E$ and let

$$X_\infty = \Delta.$$

For any subset $S$ of $E$, we write $S_\Delta$ for $S \cup \Delta$.

For any nonnegative measurable function $\Psi$ on $F_\Delta \times F_\Delta$, let us consider a random measure $\kappa(\Psi, \cdot)$ defined by

$$(3.5) \qquad \kappa(\Psi, dt) = \sum_{0 < s: R(s-) < R(s)} \Psi(X_{R(s-)}, X_{R(s)}) \varepsilon_s(dt),$$

where $\varepsilon_s$ is the point mass at $s$. By the above note, the random measure $\kappa$ may also be written as

$$\kappa(\Psi, dt) = \sum_{0 < s: s \in I} \Psi(X_s, X_{R(s)}) \varepsilon_s(dt).$$

Any function $f$ on $F$ is extended to $F_\Delta$ by setting $f(\Delta) = 0$. By this convention, $f \otimes g$ denotes the function on $F_\Delta \times F_\Delta$ defined by $\Psi(x,y) = f(x)g(y)$ for $f, g \in \mathcal{B}_+(F)$. We further let $(f \otimes I_\Delta)(x,y) = f(x)I_\Delta(y)$. Obviously, we have, for $f, g \in \mathcal{B}_+(F)$,

$$(3.6) \qquad \kappa(f \otimes g, dt) = \sum_{0 < s: R(s-) < R(s) < \infty} f(X_{R(s-)}) g(X_{R(s)}) \varepsilon_s(dt),$$

$$(3.7) \qquad \kappa(f \otimes I_\Delta, dt) = \sum_{\substack{0 < s: R(s-) < \infty \\ R(s) = \infty}} f(X_{R(s-)}) \varepsilon_s(dt).$$

For later reference, we introduce the last exit time from $F$ defined by

$$(3.8) \qquad S_F = \sup\{t > 0 : X_t \in F\}, \qquad \sup \varnothing = 0.$$

Then $s = S_F > 0$ if and only if $R(s-) < \infty$, $R(s) = \infty$ and accordingly

$$\kappa(f \otimes I_\Delta, dt) = f(X_{S_F-}) \varepsilon_{S_F}(dt).$$

LEMMA 3.2. *The random measure $\kappa(\Psi, \cdot)$ is homogeneous for any $\Psi \in \mathcal{B}_+(F_\Delta \times F_\Delta)$.*

PROOF. Since $R(s) \circ \theta_u + u = R(u+s)$, we have $X_{R(s)} \circ \theta_u = X_{R(u+s)}$ and

$$\kappa(\Psi, dt) \circ \theta_u = \sum_{u < s+u: R(u+s-) < R(u+s)} F(X_{R(u+s-)}, X_{R(u+s)}) \varepsilon_s(dt)$$

$$= \sum_{u<s:\, R(s-)<R(s)} F(X_{R(s-)}, X_{R(s)}) \varepsilon_s(dt+u)$$

$$= \kappa(\Psi, dt+u) \qquad \square$$

THEOREM 3.1. *Let $f, g \in \mathcal{B}_+(F)$. Then*

$$E^m \kappa(f \otimes g, (0,t)) = tU(f \otimes g), \qquad E^m \kappa(f \otimes I_\Delta, (0,t)) = tV(f), \qquad t > 0.$$

PROOF. For $n \geq 1$, let $D_n := \{t_{n,k} = (k-1)/2^n : k \geq 1\}$ and $I_{n,k} = [t_{n,k-1}, t_{n,k})$ for $k \geq 1$.

For $0 < s$, we observe that $R(s-) < R(s)$ and $(R(s-), R(s)) \cap D_n \neq \varnothing$ if and only if $R(s-) = L(t_{n,k}) \in I_{n,k}$, $X_{t_{n,k}} \in G$ and $R(s) = R(t_{n,k})$ for a unique $k$ depending on $n$.

Therefore, by the monotone convergence theorem and using (3.3) and (2.8), we get

$$E^m \kappa(f \otimes g, (0,t))$$

$$= E^m \sum_{0<s<t:\, R(s-)<R(s)<+\infty} f(X_{R(s-)}) g(X_{R(s)})$$

$$= \lim_n E^m \sum_{k:\, t_{n,k} \leq t} f(X_{L(t_{n,k})}) g(X_{R(t_{n,k})}) \mathbb{1}_{\{L(t_{n,k}) \in I_{n,k}, X_{t_{n,k}} \in G, R(t_{n,k}) < \infty\}}$$

$$= \lim_n \sum_{k:\, t_{n,k} \leq t} E^m f(X_{L(t_{n,k})}) g(X_{R(t_{n,k})}) \mathbb{1}_{\{L(t_{n,k}) \in I_{n,k}, X_{t_{n,k}} \in G, R(t_{n,k}) < \infty\}}$$

$$= \lim_n \sum_{k:\, t_{n,k} \leq t} \int_G E^x(T \in (0, 2^{-n}], f(X_T)) Hg(x) m_G(dx)$$

$$= \lim_n [2^n t] \int_G E^x(T \in (0, 2^{-n}], f(X_T)) Hg(x) m_G(dx)$$

$$= tU(f \otimes g),$$

where $[2^n t]$ is the largest integer dominated by $2^n t$.

In the same way, from (3.4) and (2.10) we get

$$E^m \kappa(f \otimes I_\Delta, (0,t))$$

$$= E^m \sum_{0<s<t, R(s-)<\infty, R(s)=\infty} f(X_{L(s)})$$

$$= \lim_n \sum_{k:\, t_{n,k} \leq t} E^m f(X_{L(t_{n,k})}) \mathbb{1}_{\{L(t_{n,k}) \in I_{n,k}, X_{t_{n,k}} \in G, R(t_{n,k}) = \infty\}}$$

$$= \lim_n \sum_{k:\, t_{n,k} \leq t} \int_G E^x(T \in (0, 2^{-n}], f(X_T)) q(x) m_G(dx)$$

$$= \lim_n [2^n t] \int_G E^x(T \in (0, 2^{-n}], f(X_T)) q(x) m_G(dx)$$
$$= tV(f). \qquad \Box$$

**4. Admissible measure and time changed process $Y$.** We still work with an $m$-symmetric conservative diffusion process $X$ on $E$. Let $(\mathcal{E}, \mathcal{F})$ be the associated Dirichlet form on $L^2(E; m)$.

By virtue of the transfer method (see [16, 10]), we can and shall assume without loss of generality that the Dirichlet space $(E, m, \mathcal{F}, \mathcal{E})$ is regular and $X$ is an associated strong Markov process on $E$ with continuous sample paths with infinite lifetime. By the regularity we mean that $E$ is a locally compact separable metric space, $m$ is a positive Radon measure on $E$ with full support and that $\mathcal{F} \cap C_0(E)$ is dense in $\mathcal{F}$ and in $C_0(E)$. Here $C_0(E)$ denotes the space of continuous functions on $E$ with compact support. The capacity associated with this Dirichlet form is denoted by Cap. A set $N$ with $\mathrm{Cap}(N)$ is called an $\mathcal{E}$-polar set. The phrase "$\mathcal{E}$-q.e." will mean "except for an $\mathcal{E}$-polar set."

A quasisupport of a Borel measure is a smallest quasiclosed set outside of which the measure vanishes. It is unique up to the $\mathcal{E}$-q.e. equivalence.

LEMMA 4.1. *For a closed set $F \subset E$ with $\mathrm{Cap}(F) > 0$, there exists a nontrivial positive Radon measure $\mu$ on $E$ such that $\mu$ charges no $\mathcal{E}$-polar set, $\mu(E \setminus F) = 0$ and the quasisupport of $\mu$ coincides with $F$, $\mathcal{E}$-q.e.*

PROOF. As in the preceding sections, we denote by $T$ the hitting time of $F$. Take an $m$-integrable strictly positive function $g$ on $E$ and set

$$\mu(B) = P^{g \cdot m}(X_T \in B, T < \infty), \qquad B \in \mathcal{B}(E).$$

Clearly $\mu(E \setminus F) = 0$ and $\mu$ is a nontrivial positive Radon measure charging no set of zero capacity. If a quasicontinuous function $f \in \mathcal{F}$ vanishes $\mu$-a.e., then $E^{g \cdot m}(e^{-T} f(X_T)) = 0$ which implies that the quasi-continuous function $E^{\cdot}(e^{-T} f(X_T))$ vanishes $m$-a.e. Hence $f = 0$ $\mathcal{E}$-q.e. on $F$ since the $\mathcal{E}$-q.e. point of $F$ is regular for $F$, and we can conclude on account of [11], Theorem 4.6.2, that $F$ is a quasisupport of $\mu$. $\Box$

We call a measure $\mu$ *admissible* for the closed set $F$ if it possesses the properties stated in Lemma 4.1 and its topological support $\mathrm{Supp}[\mu]$ equals $F$. The following sufficient condition for a measure $\mu$ to be admissible for $F$ can be shown in the same way as in the proof of Lemma 4.1 (see also [11], Problem 4.6.1).

LEMMA 4.2. *Let $F$ be a closed set with $\mathrm{Cap}(F) > 0$. If there exists a $\sigma$-finite measure $\mu$ with $\mathrm{Supp}[\mu] = F$ such that $F$ admits a Poisson kernel with respect to $\mu$ in the sense of Section 2, then $\mu$ is admissible for $F$.*

From now on, we consider a closed set $F$ with $\mathrm{Cap}(F) > 0$. We fix an admissible measure $\mu$ for $F$. Then $\mu$ is a smooth measure. Let $\phi(t)$ be the PCAF (positive continuous additive functional) with Revuz measure $\mu$ and let $\widetilde{F}$ be its support, namely,

$$\widetilde{F} = \{x \in E : P^x(R_\phi = 0) = 1\},$$

where

$$R_\phi = \inf\{t > 0 : \phi(t) > 0\}.$$

Then $\widetilde{F}$ is a quasisupport of $\mu$ (cf. [11], Theorem 5.1.5) and, hence, by choosing the exceptional set for $\phi$ appropriately, we may assume that

$$(4.1) \qquad \widetilde{F} \subset F, \qquad \mathrm{Cap}(F \setminus \widetilde{F}) = 0.$$

Let $\tau = (\tau_t)$ be the right-continuous inverse of $\phi$:

$$(4.2) \qquad \tau_t = \inf\{s : \phi(s) > t\}, \qquad \inf \varnothing = \infty.$$

We set

$$(4.3) \qquad Y_t = X_{\tau_t}, \qquad t < \check{\zeta}, \text{ where } \check{\zeta} = \phi(\infty).$$

Then $Y = (Y_t, \check{\zeta}, P^x)_{x \in \widetilde{F}}$ is a right process on the state space $\widetilde{F}$ with lifetime $\check{\zeta}$, which is called a time change of $X$ or the time-changed process (cf. [18]). We add a cemetery $\Delta$ to $\widetilde{F}$ and define

$$Y_t = \Delta, \qquad t \geq \check{\zeta},$$

so that the time-changed process $Y$ is a right process on $\widetilde{F}_\Delta = \widetilde{F} \cup \Delta$. We also note that

$$(4.4) \qquad Y_{t-} = X_{\tau_{t-}} \in F, \qquad t \leq \check{\zeta},$$

owing to the continuity of the sample path of $X$.

In general, the process $Y = (Y_t, \check{\zeta}, P^x)_{x \in \widetilde{F}}$ is not a Hunt process. It could happen that $Y_{t-} \in F \setminus \widetilde{F}$ and $Y$ may not be quasi-left continuous either. By making use of a general reduction theorem formulated in Section 8, however, we can show that the restriction of $Y$ to the outside of a suitable exceptional set is actually a Hunt process.

To this end, we recall some basic facts about the time-changed process $Y$ on $\widetilde{F}$ shown in [11], Theorem 6.2.1. Process $Y$ is $\mu$-symmetric and the associated Dirichlet form [denoted by $(\check{\mathcal{E}}, \check{\mathcal{F}})$] on $L^2(F, \mu)$ is regular. Furthermore, $Y$ is properly associated with $(\check{\mathcal{E}}, \check{\mathcal{F}})$ in the sense that $\check{p}_t u$ is an $\check{\mathcal{E}}$-quasicontinuous version of $\check{T}_t u$ for any $u \in L^2(F; \mu)$, where $\check{p}_t$ (resp. $\check{T}_t$) denotes the transition function of $Y$ [resp. the $L^2$ semigroup associated with $(\check{\mathcal{E}}, \check{\mathcal{F}})$].

It is also clear from the preceding definition of the path $Y_t$ that the left limit $Y_{t-}$ exists in $F_\Delta$ for all $t > 0$. Hence all the conditions in Theorem 8.1 are satisfied by the time-changed process $Y$ and we are led to the next theorem for $Y$. The capacity on $F$ associated with $(\check{\mathcal{E}}, \check{\mathcal{F}})$ is denoted by $\check{\mathrm{Cap}}$. A set $N \subset F$ with $\check{\mathrm{Cap}}(N) = 0$ is called an $\check{\mathcal{E}}$-polar set.

THEOREM 4.1. *There exists a Borel subset $\check{F}$ of $\widetilde{F}$ such that*

(4.5) $$F \setminus \check{F} \text{ is } \check{\mathcal{E}}\text{-polar and } \mathcal{E}\text{-polar},$$

*$\check{F}$ is $Y$-invariant and the restriction $Y|_{\check{F}}$ of the time-changed process $Y$ to $\check{F}$ is a Hunt process properly associated with $\check{\mathcal{E}}$.*

By the general theorem, Theorem 8.1, we know only that the set $F \setminus \check{F}$ is $\check{\mathcal{E}}$-polar. However, then $\widetilde{F} \setminus \check{F}$ is $\mathcal{E}$-polar by virtue of [11], Lemma 6.2.5. Hence
$$F \setminus \check{F} = (F \setminus \widetilde{F}) + (\widetilde{F} \setminus \check{F})$$
is $\mathcal{E}$-polar as well in view of (4.1).

Finally we notice that the Dirichlet form $(\check{\mathcal{E}}, \check{\mathcal{F}})$ admits the following description. Denote by $\mathcal{F}_e$ the extended Dirichlet space of $\mathcal{F}$ and take any $u \in \mathcal{F}_e$ to be $\mathcal{E}$-quasicontinuous. Then, due to (4.1) and [11], Theorem 6.2.1,

(4.6)
$$\check{\mathcal{F}} = \{f \in L^2(F; \mu) : f = u \ \mu\text{-a.e. on } F \text{ for some } u \in \mathcal{F}_e\},$$
$$\check{\mathcal{E}}(f, f) = \mathcal{E}(Hu, Hu), \qquad f \in \check{\mathcal{F}}, \ f = u \ \mu\text{-a.e. on } F, \ u \in \mathcal{F}_e,$$

where $Hu$ is defined by
$$Hu(x) = E^x(u(X_T); T < \infty), \qquad x \in E.$$

**5. Identification of jumping and killing measures of $Y$ with $U$ and $V$.** For simplicity, the restriction of the time-changed process $Y$ to the set $\check{F}$ of Theorem 4.1 is again denoted by $Y$. Then $Y$ is a Hunt process on $\check{F} \cup \Delta$ properly associated with the regular Dirichlet form $(\check{\mathcal{E}}, \check{\mathcal{F}})$ on $L^2(F; \mu)$ and $F \setminus \check{F}$ is not only $\check{\mathcal{E}}$-polar, but also $\mathcal{E}$-polar.

Since the Dirichlet form $(\check{\mathcal{E}}, \check{\mathcal{F}})$ on $L^2(F; \mu)$ is regular, it admits the Beurling–Deny decomposition; for any $\check{\mathcal{E}}$-quasicontinuous functions $f, g \in \check{\mathcal{F}}$,

(5.1)
$$\check{\mathcal{E}}(f, g) = \check{\mathcal{E}}^{(c)}(f, g) + \int_{F \times F \setminus d} (f(x) - f(y))(g(x) - g(y)) J(dx, dy)$$
$$+ \int_F f(x) g(x) k(dx),$$

where $\check{\mathcal{E}}^{(c)}$ is a symmetric form with a strong local property, $J$ is a symmetric positive Radon measure on $F \times F$ off the diagonal $d$ and $k$ is a positive Radon measure on $F$. Measures $J$ and $k$ are called the *jumping measure* and the *killing measure* for the Dirichlet form $(\check{\mathcal{E}}, \check{\mathcal{F}})$, respectively.

Since $Y$ is a Hunt process properly associated with $(\check{\mathcal{E}}, \check{\mathcal{F}})$, we can use directly the general result of [11], Section 5.3, to describe $J$ and $k$ in terms of the Lévy system of $Y$. Let $(N(x, dy), \psi)$ be a Lévy system of $Y$. More precisely $N(x, dy)$ is a kernel on $(\check{F}_\Delta, \mathcal{B}(\check{F}_\Delta))$ with $N(x, \{x\}) = 0$, $x \in \check{F}$, and $\psi = \psi(t)$ is a PCAF of $Y$ such that, for any $\Psi \in \mathcal{B}^+(\check{F}_\Delta \times \check{F}_\Delta)$ vanishing on the diagonal,

$$
\begin{aligned}
E^x & \left( \sum_{s \leq t} \Psi(Y_{s-}, Y_s) \right) \\
& = E^x \left( \int_0^t \int_{\check{F}_\Delta} N(Y_s, dy) \Psi(Y_s, y) \, d\psi(s) \right), \qquad x \in \check{F}.
\end{aligned}
\tag{5.2}
$$

Let $\nu$ be the Revuz measure of $\psi$ with respect to $Y$. Then, by [11], Theorem 5.3.1,

$$
J(dx, dy) = \tfrac{1}{2} N(x, dy) \nu(dx), \qquad k(dx) = N(x, \Delta) \nu(dx). \tag{5.3}
$$

By the Revuz correspondence, we have, for any $\Psi \in \mathcal{B}^+(\check{F}_\Delta \times F_\Delta)$ vanishing on the diagonal,

$$
\begin{aligned}
\int_{\check{F} \times \check{F} \setminus d} \Psi(x, y) J(dx, dy) &= \lim_{t \downarrow 0} \frac{1}{2t} E^\mu \sum_{0 < s \leq t} \Psi(Y_{s-}, Y_s) I_{\check{F}}(Y_s) \\
&= \lim_{\alpha \to \infty} \frac{\alpha}{2} E^\mu \sum_{0 < t < \infty} e^{-\alpha t} \Psi(Y_{t-}, Y_t) I_{\check{F}}(Y_t),
\end{aligned}
\tag{5.4}
$$

$$
\begin{aligned}
\int_{\check{F}} \Psi(x, \Delta) k(dx) &= \lim_{t \downarrow 0} \frac{1}{t} E^\mu \sum_{0 < s \leq t} \Psi(Y_{s-}, \Delta) I_\Delta(Y_s) \\
&= \lim_{\alpha \to \infty} \alpha E^\mu \sum_{0 < t < \infty} e^{-\alpha t} \Psi(Y_{t-}, \Delta) I_\Delta(Y_t).
\end{aligned}
\tag{5.5}
$$

THEOREM 5.1. *We have*

(5.6) $\quad J \geq \tfrac{1}{2} U \qquad$ *on* $\check{F} \times \check{F} \setminus d$, $\qquad U = 0 \qquad$ *on* $(F \times F \setminus d) \setminus (\check{F} \times \check{F})$,

(5.7) $\quad k \geq V \qquad$ *on* $\check{F}$, $\qquad\qquad\qquad V = 0 \qquad$ *on* $F \setminus \check{F}$.

*Furthermore, if there exists a sequence of finite $X$-excessive measures $m_n$ increasing to $m$, then $J = \tfrac{1}{2} U$ and $k = V$.*

PROOF. It is known that $R_\phi = \widetilde{T}$, where $\widetilde{T}$ is the hitting time of the support $\widetilde{F}$ of $\phi$. Hence

$$\tau_{\phi(t)} = \inf\{s : \phi(s) > \phi(t)\} = \inf\{s > t : \phi(s-t) \circ \theta_t > 0\} = \widetilde{T} \circ \theta_t + t.$$

Since $F \setminus \widetilde{F}$ is $\mathcal{E}$-polar, we have

$$P^x(\widetilde{T} = T) = 1, \qquad \mathcal{E}\text{-q.e. } x \in E,$$

and hence

(5.8) $\qquad \tau_{\phi(t)} = R(t) \qquad \forall t > 0, \ P^x\text{-a.e. for } \mathcal{E}\text{-q.e. } x \in E.$

For any $\Psi \in \mathcal{B}^+(F \times F)$ vanishing on the diagonal, we have from (5.4), (4.4) and (5.8),

$$2 \int_{\check{F} \times \check{F} \setminus d} \Psi(x,y) J(dx,dy) = \lim_{\alpha \to \infty} \alpha E^\mu \sum_{0 < t < \infty} e^{-\alpha t} \Psi(Y_{t-}, Y_t) I_{\check{F}}(Y_t)$$

$$= \lim_{\alpha \to \infty} \alpha E^\mu \sum_{0 < t < \infty} e^{-\alpha t} \Psi(X_{\tau_{t-}}, X_{\tau_t}) I_{\widetilde{F}}(X_{\tau_t})$$

$$= \lim_{\alpha \to \infty} \alpha E^\mu(\Sigma_\alpha),$$

where

$$\Sigma_\alpha = \sum_{\substack{0 < t < \infty \\ R(t) < \infty}} e^{-\alpha \phi(t)} \Psi(X_{R(t-)}, X_{R(t)}).$$

Since $\mu$ is the Revuz measure of $\phi$ with respect to the conservative $m$-symmetric process $X$, we have from [11], Theorem 5.1.3, and [18], (32.6), that

$$\alpha E^\mu(\Sigma_\alpha) = \alpha \frac{1}{s} E^m \left( \int_0^s E^{X_u}(\Sigma_\alpha) \, d\phi(u) \right) = \alpha \frac{1}{s} E^m \int_0^s \Sigma_\alpha \circ \theta_u \, d\phi(u)$$

$$= \alpha \frac{1}{s} E^m \int_0^s \sum_{\substack{0 < t < \infty \\ R(t+u) < \infty}} e^{-\alpha(\phi(t+u) - \phi(u))} \Psi(X_{R(t+u-)}, X_{R(t+u)}) \, d\phi(u)$$

$$= \alpha \frac{1}{s} E^m \int_0^s e^{\alpha \phi(u)} \, d\phi(u) \sum_{\substack{u < t < \infty \\ R(t) < \infty}} e^{-\alpha \phi(t)} \Psi(X_{R(t-)}, X_{R(t)})$$

$$= \frac{1}{s} E^m \left( \sum_{R(t) < \infty} e^{-\alpha \phi(t)} \Psi(X_{R(t-)}, X_{R(t)}) \int_0^s I_{\{t > u\}}(u) \, de^{\alpha \phi(u)} \right),$$

$$= \frac{1}{s} E^m \left( \sum_{R(t) < \infty} e^{-\alpha \phi(t)} \Psi(X_{R(t-)}, X_{R(t)}) \cdot (e^{\alpha \phi(s \wedge t)} - 1) \right)$$

$$= \frac{1}{s}(I_\alpha^- + I_\alpha^+),$$

where

$$I_\alpha^- := E^m\left[\sum_{\substack{0<t\leq s \\ R(t)<\infty}} (1 - e^{-\alpha\phi(t)})\Psi(X_{R(t-)}, X_{R(t)})\right],$$

$$I_\alpha^+ := E^m\left[\sum_{\substack{s<t<\infty \\ R(t)<\infty}} e^{-\alpha\phi(t)}\Psi(X_{R(t-)}, X_{R(t)})(e^{\alpha\phi(s)} - 1)\right],$$

and they both are nonnegative. In other words, we have

$$2\int_{\check{F}\times\check{F}\setminus d} \Psi\, dJ = \lim_{\alpha\to\infty} \frac{1}{s}(I_\alpha^- + I_\alpha^+),$$

as an increasing limit. By the monotone convergence theorem, we arrive at the inequality

(5.9)
$$2\int_{\check{F}\times\check{F}\setminus d} \Psi dJ \geq \frac{1}{s}\lim_{\alpha\to\infty} I_\alpha^-$$
$$= \frac{1}{s}E^m\left[\sum_{\substack{0<t\leq s \\ R(t)<\infty}} \Psi(X_{R(t-)}, X_{R(t)})\right].$$

Here we note that we can insert condition $\phi(t) > 0$ in the summand of $I_\alpha^-$, because this is true for $t > T$ while $\phi(T) = 0$ but the path can not be the left end point of an excursion at time $T$ because q.e. point of $F$ is regular for itself. By Theorem 3.1, (5.6) holds. The proof of (5.7) is similar.

Assume now the existence of finite $X$-excessive measures $\{m_n\}$ increasing to $m$. This assumption is trivially satisfied when $m(E) < \infty$. Choose any nonnegative function $\Psi$ vanishing on $d$ such that $\int_{\check{F}\times\check{F}\setminus d} \Psi\, dJ < \infty$. We let, for any $u \geq 0$,

$$\Sigma_{\alpha,u} := \sum_{\substack{u<t<\infty \\ R(t)<\infty}} e^{-\alpha\phi(t)}\Psi(X_{R(t-)}, X_{R(t)}),$$

so that $I_\alpha^+ = E^m(\Sigma_{\alpha,s}(e^{\alpha\phi(s)} - 1))$. It can be easily verified that $\Sigma_{\alpha,s} \cdot e^{\alpha\phi(s)} = \Sigma_{\alpha,0} \circ \theta_s$.

Take a truncation function $\chi_N(x) = x \wedge N$, $x \in \mathbb{R}$, and set

$$I_{\alpha,n,N}^+ = E^{m_n}(\chi_N(\Sigma_{\alpha,s} \cdot e^{\alpha\phi(s)}) - \chi_N(\Sigma_{\alpha,s})),$$

which then increases to $I_\alpha^+$ when we let $N \uparrow \infty$ and then $n \uparrow \infty$. Since $m_n$ is a finite excessive measure, we have

$$I_{\alpha,n,N}^+ = E^{m_n}(E^{X_s}(\chi_N(\Sigma_{\alpha,0}))) - E^{m_n}(\chi_N(\Sigma_{\alpha,s}))$$

$$\leq E^{m_n}(\chi_N(\Sigma_{\alpha,0})) - E^{m_n}(\chi_N(\Sigma_{\alpha,s}))$$
$$= E^{m_n}(\chi_N(\Sigma_{\alpha,0}) - \chi_N(\Sigma_{\alpha,s}))$$
$$= E^{m_n}(\Sigma_{\alpha,0} - \Sigma_{\alpha,s}; N > \Sigma_{\alpha,0}) + E^{m_n}(N - \Sigma_{\alpha,0}; \Sigma_{\alpha,0} \geq N > \Sigma_{\alpha,s})$$
$$\leq E^{m_n}(\Sigma_{\alpha,0} - \Sigma_{\alpha,s})$$
$$\leq E^m \left( \sum_{\substack{0 < t \leq s \\ R(t) < \infty,\ \phi(t) > 0}} e^{-\alpha\phi(s)} \Psi(X_{R(t-)}, X_{R(t)}) \right).$$

The last expectation in the above is finite in view of (5.9) and (5.10).

It follows that

$$I_\alpha^- + I_{\alpha,n,N}^+ \leq E^m \left( \sum_{\substack{0 < t \leq s \\ R(t) < \infty}} \Psi(X_{R(t-)}, X_{R(t)}) \right).$$

Therefore we have

(5.10)
$$2 \int_{\check{F} \times \check{F} \setminus d} \Psi\, dJ = \lim_{\alpha \to \infty} \frac{1}{s}(I_\alpha^- + \lim_n \lim_N I_{\alpha,n,N}^+)$$
$$\leq \frac{1}{s} E^m \left( \sum_{\substack{0 < t \leq s \\ R(t) < \infty}} \Psi(X_{R(t-)}, X_{R(t)}) \right).$$

Inequalities (5.9) and (5.10) give an equality

$$2 \int_{\check{F} \times \check{F} \setminus d} \Psi\, dJ = \frac{1}{s} E^m \left( \sum_{\substack{0 < t \leq s \\ R(t) < \infty}} \Psi(X_{R(t-)}, X_{R(t)}) \right).$$

By substituting $\Psi = f \otimes g$ for any $f, g \in C_0(F)$ with disjoint support in the above equality, we get the desired identity by virtue of Theorem 3.1. The proof for $k = V$ is similar. □

COROLLARY 5.1. *Suppose that $m$ is finite and the hitting measure has the Poisson kernel $K(x, \xi), x \in G, \xi \in F$ with respect to a $\sigma$-finite measure $\mu$ with $\mathrm{Supp}[\mu] = F$. Then $\mu$ is admissible and the associated time-changed process $Y$ (with a possible q.e. modification of its state space $\check{F}$) has as its Lévy system*

(5.11)
$$(U(\xi, \eta)\mu(d\eta), t),$$

*where $U$ is the Feller kernel defined by (2.13) in terms of $K$ and $t$ denotes the nonrandom PCAF $\psi(t) = t$ of $Y$.*

PROOF. Measure $\mu$ is admissible by Lemma 4.2. By Theorem 5.1,

$$J(d\xi, d\eta) = \tfrac{1}{2} U(\xi, \eta) \mu(d\xi) \mu(d\eta) \qquad \text{on } \check{F} \times \check{F} \setminus d$$

and $\mu$ is the Revuz measure of the PCAF $t$ of $Y$. Hence, it suffices to show that the value of the right-hand side of (5.2) depends only on the function $\Psi$ and the jumping measure $J$ for q.e. $x$ and that it does not depend on the special choice of $N$ and $\nu$ that express $J$ as in (5.3). This can been readily seen from known formulae ([11], (5.1.12) and (5.1.14)) on the Revuz correspondence of the PCAF and the smooth measure. □

**6. Trace Dirichlet form and Douglas integral with Feller measure.** In the preceding two sections, we proved the following: Let $(E, m, \mathcal{F}, \mathcal{E})$ be a regular Dirichlet space and let $X$ be an associated conservative diffusion process on $E$. Any function in the extended Dirichlet space $\mathcal{F}_e$ is taken to be $\mathcal{E}$-quasicontinuous. Let $F$ be a closed subset of $X$ with $\mathrm{Cap}(F) > 0$, let $\mu$ be an admissible measure for $F$ with $\mathrm{Supp}[\mu] = F$, let $\phi$ be a PCAF of $X$ with Revuz measure $\mu$ and let $Y$ be a time-changed process of $X$ by means of $\phi$. Process $Y$ is $\mu$-symmetric and its Dirichlet form on $L^2(F; \mu)$ is denoted by $(\check{\mathcal{E}}, \check{\mathcal{F}})$ which is also called the *trace Dirichlet form* of $\mathcal{E}$ on $F$. In fact, in view of (4.6) and [11], Lemma 6.2.5, we have

(6.1)
$$\check{\mathcal{F}} = \mathcal{F}_e|_F \cap L^2(F; \mu),$$
$$\check{\mathcal{E}}(f, f) = \mathcal{E}(Hu, Hu), \qquad f = u|_F, \ u \in \mathcal{F}_e.$$

Furthermore we have obtained in Theorem 5.1 that, for any $f, g \in \check{\mathcal{F}}$,

(6.2)
$$\check{\mathcal{E}}(f, g) \geq \check{\mathcal{E}}^{(c)}(f, g) + \tfrac{1}{2} \int_{F \times F \setminus d} (f(x) - f(y))^2 (g(x) - g(y)) U(dx, dy)$$
$$+ \int_F f(x)^2 g(x) V(dx),$$

the representation of the trace Dirichlet form $\check{\mathcal{E}}$ in terms of the Feller measure $U$ and the supplementary Feller measure $V$ introduced in Section 2. In particular, the first integral on the right-hand side is called the *Douglas integral* with the Feller measure $U$.

The Feller measure $U$ and the supplementary Feller measure $V$ are completely determined by the absorbed (minimal) process $X_G$ of $X$ on the set $G = X \setminus F$, while the local term $\check{\mathcal{E}}^{(c)}$ in the above decomposition is determined by the behavior of $X$ on the set $F$. On the other hand, the value of the Dirichlet form $\mathcal{E}(u, u)$ for $u \in \mathcal{F}_e$ is known to be equal to half of the total mass of the energy measure $\mu_{\langle u \rangle}$ of $u$. Therefore we may expect that

(6.3) $\qquad \check{\mathcal{E}}^{(c)}(f, f) = \tfrac{1}{2} \mu_{\langle Hu \rangle}(F), \qquad f = u|_F, \ u \in \mathcal{F}_e.$

We do not prove this, but more specifically, we show in this section that if $\mu_{\langle u \rangle}(F)$ vanishes for any $u \in \mathcal{F}_e$, then the trace Dirichlet form equals the Douglas integral with the Feller measure under the assumption that $m(G)$ is finite.

To this end, we first show the domination of the trace Dirichlet form by the Douglas integral under the setting that $(E, m, \mathcal{F}, \mathcal{E})$ is a regular irreducible Dirichlet space and $X$ is an associated Hunt process on $E$. We do not assume that $X$ is of continuous sample paths, but we assume that $X$ is conservative. We further assume that

$$(6.4) \qquad m(G) < \infty, \qquad \text{Cap}(F) > 0.$$

We note that (6.4) and the irreducibility of $\mathcal{E}$ imply that

$$(6.5) \qquad P^x(T < \infty) = 1, \qquad \text{q.e. } x \in G,$$

because then $P^x(T < \infty) > 0$ for q.e. $x \in E$ by [11], Theorem 4.6.6, and Lemma 2.3 applies.

For any $u \in \mathcal{F}_{e,b} = \mathcal{F}_e \cap L^\infty(E; m)$, its energy measure $\mu_{\langle u \rangle}$ is defined by

$$(6.6) \qquad \int_E f(x) \mu_{\langle u \rangle}(dx) = 2\mathcal{E}(uf, u) - \mathcal{E}(u^2, f), \qquad f \in \mathcal{F} \cap C_0(E).$$

The energy measure uniquely extends to any $u \in \mathcal{F}_e$ and it holds that

$$(6.7) \qquad \mathcal{E}(u, u) = \tfrac{1}{2} \mu_{\langle u \rangle}(E), \qquad u \in \mathcal{F}_e.$$

Let

$$\mathcal{F}^0 = \{u \in \mathcal{F} : u = 0 \text{ q.e. on } F\}.$$

Then $(\mathcal{F}^0, \mathcal{E})$ is a regular Dirichlet space on $L^2(G; m)$ which is associated with the absorbed process $X_G$ ([11], Theorem 4.4.3). Recall that $R^0_\alpha$ denotes the resolvent operator for $X_G$. Since

$$R^0_1 1(x) = 1 - E^x(e^{-T}) < 1, \qquad \text{q.e. } x \in G,$$

by (6.5), we see that $(\mathcal{F}^0, \mathcal{E})$ is transient by virtue of [11], Lemma 1.6.5, and moreover, the extended Dirichlet space $\mathcal{F}^0_e$ of $\mathcal{F}^0$ admits the expression

$$\mathcal{F}^0_e = \{u \in \mathcal{F}_e : u = 0 \text{ q.e. on } F\}$$

due to [11], Theorem 4.4.4.

Denote by $S^{(0)}_0(G)$ the space of positive Radon measures of finite 0-order energy integral with respect to $(\mathcal{F}^0_e, \mathcal{E})$. If $\nu \in S^{(0)}_0(G)$, then there exists a unique $R^0 \nu \in \mathcal{F}^0_e$ called the 0-order potential of $\nu$ such that

$$(6.8) \qquad \mathcal{E}(R^0 \nu, v) = \int_G v \, d\nu, \qquad v \in \mathcal{F} \cap C_0(G).$$

Equation (6.7) extends to any quasicontinuous function $v \in \mathcal{F}_e^0$.

We write $(f,g)_G = \int_G fg\,dm$. We know from [11], Theorem 1.5.4, that if a nonnegative measurable function $f$ on $G$ satisfies that $(f, R_{0+}^0 f)_G < \infty$, then $R_{0+}^0 f \in \mathcal{F}_e^0$ and

$$\mathcal{E}(R_{0+}^0 f, v) = (f, v)_G, \qquad v \in \mathcal{F}_e^0. \tag{6.9}$$

We know further from [11], Theorem 4.6.5, that $Hu \in \mathcal{F}_e$ for any $u \in \mathcal{F}_e$ and

$$\mathcal{E}(Hu, v) = 0 \qquad \forall v \in \mathcal{F}_e^0. \tag{6.10}$$

We prepare a lemma which generalizes the methods of computing the Dirichlet norms of classical harmonic functions employed in [4] and [8].

LEMMA 6.1. *For any $u \in \mathcal{F}_{e,b}$, let*

$$w = H(u^2) - (Hu)^2 \ (\in \mathcal{F}_{e,b}).$$

*Then*

$$w \in \mathcal{F}_{e,b}^0 \quad \text{and} \quad w = R^0 \nu \qquad \text{for } \nu = \mu_{\langle Hu \rangle}|_G. \tag{6.11}$$

*Furthermore,*

$$\mu_{\langle Hu \rangle}(G) = \lim_{\alpha \to \infty} \alpha(H^\alpha 1, w)_G. \tag{6.12}$$

PROOF. Since $\mathcal{F}_{e,b}$ is an algebra and $w = 0$ q.e. on $F$, we have that $w \in \mathcal{F}_{e,b}^0$. From (6.6) and (6.10), we have, for any $f \in \mathcal{F} \cap C_0(D)(\subset \mathcal{F}^0)$,

$$\mathcal{E}(w, f) = -\mathcal{E}((Hu)^2, f)$$
$$= 2\mathcal{E}(Hu \cdot f, Hu) - \mathcal{E}((Hu)^2, f) = \int_G f\,d\mu_{\langle Hu \rangle},$$

arriving at (6.11).

Since $H1 = 1$ q.e. on $G$ by (6.5), we have

$$1 - H^\alpha 1 = \alpha R_{0+}^0 H^\alpha 1$$

and hence $(H^\alpha 1, R_{0+}^0 H^\alpha 1)_G < \infty$. Accordingly, for any $\nu \in S_0^{(0)}(G)$, from (6.8) and (6.9) we get

$$\alpha(H^\alpha 1, R^0\nu)_G = \mathcal{E}(\alpha R_{0+}^0 H^\alpha 1, R^0\nu) = \langle \alpha R_{0+}^0 H^\alpha 1, \nu \rangle = \langle 1 - H^\alpha 1, \nu \rangle \uparrow \nu(G).$$

Hence (6.12) follows from (6.11). □

THEOREM 6.1. *For any $u \in \mathcal{F}_e$,*

$$\mu_{\langle Hu \rangle}(G) \leq \int_{F \times F} (u(\xi) - u(\eta))^2 U(d\xi, d\eta). \tag{6.13}$$

PROOF. This follows from (6.12) and the identity in [8], (15),

$$\alpha(H^\alpha 1, w)_G + \alpha \int_{F \times G} (Hu(x) - u(\xi))^2 H^\alpha(x, d\xi) m(dx)$$
$$(6.14) \qquad = \int_{F \times F} (u(\xi) - u(\eta))^2 U_\alpha(d\xi, d\eta),$$

which can be easily verified. □

Theorem 6.1 combined with (6.7) leads to the following.

COROLLARY 6.1. *Suppose that*

$$(6.15) \qquad \mu_{\langle u \rangle}(F) = 0 \qquad \forall u \in \mathcal{F}_e.$$

*Then, for any* $u \in \mathcal{F}_e$,

$$(6.16) \qquad \mathcal{E}(Hu, Hu) \leq \tfrac{1}{2} \int_{F \times F} (u(\xi) - u(\eta))^2 U(d\xi, d\eta).$$

We emphasize that condition (6.15) is satisfied if the energy measure of $u$ is absolutely continuous with respect to $m$, that is, *a carré du champ operator* $\Gamma(u, u)$ *exists for any* $u \in \mathcal{F}$ *and* $m(F) = 0$.

We can now state the main theorem of this section.

THEOREM 6.2. *Let* $(E, m, \mathcal{F}, \mathcal{E})$ *be a regular irreducible Dirichlet space whose associated Markov process on* $E$ *is a conservative diffusion. For a closed set* $F \subset E$ *and its complement* $G$, *we assume condition* (6.4). *We further assume condition* (6.15) *for the energy measures associated with* $\mathcal{E}$. *Then, for any* $u \in L^2(F, \mu) \cap \mathcal{F}_e$,

$$(6.17) \qquad \mathcal{E}(Hu, Hu) = \tfrac{1}{2} \int_{F \times F} (u(\xi) - u(\eta))^2 U(d\xi, d\eta).$$

PROOF. By (6.1) and (6.2), we already have the converse inequality to (6.16).

□

We may view Theorem 6.2 from a quite different angle. The Dirichlet form $(\mathcal{F}, \mathcal{E})$ on $L^2(E; m)$ is in a sense an extension of the absorbed Dirichlet space $(\mathcal{F}^0, \mathcal{E})$ on $L^2(G; m)$. What kind of extension are we dealing with under condition (6.15)? This question can be answered in terms of the notion of the reflected Dirichlet space initially formulated by Silverstein [19, 20] and finally by Chen [3].

We continue to consider a regular irreducible Dirichlet space $(E, m, \mathcal{F}, \mathcal{E})$ associated with a conservative diffusion $X$ on $E$ and we assume condition (6.4) for a closed set $F \subset E$ and its complement $G$.

Let $(\mathcal{F}_a^{\mathrm{ref}}, \mathcal{E}^{\mathrm{ref}})$ be the $L^2$ *reflected Dirichlet space* (in the sense of [3]) relative to the regular Dirichlet space $(\mathcal{F}^0, \mathcal{E})$ on $L^2(G; m)$ associated with the absorbed process $X_G$.

THEOREM 6.3. *Condition* (6.15) *is equivalent to*

$$(6.18) \qquad \mathcal{F}|_G \subset \mathcal{F}_a^{\mathrm{ref}}, \qquad \mathcal{E}(u, v) = \mathcal{E}^{\mathrm{ref}}(u|_G, v|_G), \qquad u, v \in \mathcal{F}.$$

PROOF. By (6.10) and the preceding description of the space $\mathcal{F}_e^0$, we have, for any $u \in \mathcal{F}_e$,

$$(6.19) \quad u_0 = u - Hu \in \mathcal{F}_e^0, \qquad \mathcal{E}(u, u) = \mathcal{E}(u_0, u_0) + \mathcal{E}(Hu, Hu).$$

We further know from (6.7) that condition (6.15) is equivalent to

$$(6.20) \qquad \mathcal{E}(Hu, Hu) = \tfrac{1}{2} \mu_{\langle Hu \rangle}(G) \qquad \forall u \in \mathcal{F}_e.$$

Let $G_k$ be relatively compact open sets increasing to $G$ and let $L_k$ be the equilibrium measures of the 0-order equilibrium potentials $e_k$ for the sets $G_k$ relative to the extended Dirichlet space $(\mathcal{F}_e^0, \mathcal{E})$:

$$e_k \in \mathcal{F}_e^0, \ 0 \leq e_k \leq 1, \ e_k = 1 \text{ on } G_k, \qquad \mathcal{E}(e_k, v) = \langle v, L_k \rangle_G, \ v \in \mathcal{F}_e^0.$$

We then have

$$(6.21) \qquad \mu_{\langle Hu \rangle}(G) = \lim_{k \to \infty} \langle H(u^2) - (Hu)^2, L_k \rangle_G, \qquad u \in \mathcal{F}_{e,b}.$$

In fact, using the notation in Lemma 6.1 we see that

$$\langle w, L_k \rangle_G = \mathcal{E}(w, e_k) = \langle e_k, \nu \rangle_G,$$

which tends as $k \to \infty$ to $\nu(G) = \mu_{\langle Hu \rangle}(G)$. By comparing the combination of (6.19), (6.20) and (6.21) with Definition 3.1 in [3] of the $L^2$ reflected Dirichlet space, we get the equivalence of (6.15) and (6.18). □

Takeda ([22], Theorem 3.3) showed that the $L^2$ reflected Dirichlet space is the maximum Silverstein extension of $(\mathcal{F}^0, \mathcal{E})$ in a specific semiorder. When $(\mathcal{F}^0, \mathcal{E})$ is the Dirichlet space of the absorbing Brownian motion on an arbitrary bounded domain $D$, $\mathcal{F}_a^{\mathrm{ref}}$ equals $H^1(D)$, which has been described [9] in terms of the Feller kernel on the Martin boundary (see also the next section). In view of Theorem 6.3, we thus see that the Dirichlet space $(\mathcal{F}, \mathcal{E})$ satisfying condition (6.15) corresponds to a member of the class $\mathbf{G}_1$ of [9], Section 8, in this special case.

**7. Application to reflecting Brownian motion on a Lipschitz domain.** Let $D$ be a bounded Lipschitz domain of $\mathbb{R}^d$ with $d \geq 2$ and let $\overline{D} = D \cup \partial D$ be its closure. As is well known [2, 15], $\partial D$ (resp. $\overline{D}$) can then be identified with the Martin boundary $M$ of $D$ (resp. the Martin space $D \cup M$) and $M$ consists only of the minimal boundary points. In what follows, we regard the relative boundary $\partial D$ also as the Martin boundary of $D$ under this identification.

Denote a Martin kernel by $K(x,\xi)$, $x \in D, \xi \in \partial D$. By the Martin representation theorem [5], any positive harmonic function $h$ on $D$ can be expressed as the integral of the Martin kernel against a unique positive Radon measure on $\partial D$ called the *Martin representing measure of $h$ corresponding to $K$*. We let $\mu$ be the Martin representing measure of the constant harmonic function 1 corresponding to $K$:

$$(7.1) \qquad 1 = \int_{\partial D} K(x,\xi)\mu(d\xi), \qquad x \in D.$$

We now consider the space

$$(7.2) \quad \mathcal{F} = H^1(D), \qquad \mathcal{E}(u,v) = \tfrac{1}{2}\int_D \nabla u \cdot \nabla v\, dx, \qquad u, v \in H^1(D),$$

which can be regarded as a strongly local regular Dirichlet space on $L^2(\overline{D})$ [rather than $L^2(D)$] and hence there exists an associated conservative diffusion process $X = (X_t, P^x)_{x \in \overline{D}}$ on $\overline{D}$ uniquely up to the q.e. equivalence ([11], Example 4.5.3). We fix such a process $X$ and call it a *reflecting Brownian motion* on $\overline{D}$.

Let $T$ be the hitting time of $\partial D$ on $X$ and let $H(x,\cdot)$ be the hitting distribution of $X$ on $\partial D$:

$$H(x, B) = P^x(X_T \in B, T < \infty), \qquad x \in D,\ B \in \mathcal{B}(\partial D).$$

LEMMA 7.1. *The hitting distribution $H(x,\cdot)$ of $X$ and the measure $\mu$ in (7.1) are related by*

$$H(x,B) = \int_B K(x,\xi)\mu(d\xi) \qquad \forall B \in \mathcal{B}(\partial D)\ \textit{for q.e. } x \in D.$$

PROOF. Let $X_D$ be the absorbed process of $X$ obtained by killing $X$ at time $T$. Thus, $X_D = (X_t, P^x, \zeta^D)$ with lifetime $\zeta^D$ given by

$$(7.3) \qquad\qquad\qquad \zeta^D = T.$$

By virtue of [11], Theorem 4.4.3, $X_D$ is associated with the part of the Dirichlet form (7.2) on the open set $D$, namely

$$(7.4)\quad \mathcal{F}_D = H^1_0(D), \qquad \mathcal{E}_D(u,v) = \tfrac{1}{2}\int_D \nabla u \cdot \nabla v\, dx, \qquad u, v \in H^1_0(D).$$

Since the absorbing Brownian motion on $D$ (the standard Brownian motion on $\mathbb{R}^d$ killed upon leaving the set $D$) is also associated with the Dirichlet form (7.4) ([11], Example 4.4.1), we see that $X_D = (X_t, P^x, \zeta^D)$ coincides in law with the absorbing Brownian motion on $D$ for q.e. starting point $x \in D$.

According to Doob's description of the structure of Brownian motion on the Martin space ([5], page 727), we therefore have that

$$(7.5)\ P^x(X_{\zeta^D-} \in B) = \int_B K(x,\xi)\mu(d\xi) \qquad \forall B \in \mathcal{B}(\partial D) \text{ for q.e. } x \in D.$$

The lemma follows from (7.3) and (7.5). $\square$

Since the Martin kernel $K(x,\xi)$ is harmonic in $x \in D$, it is excessive with respect to the absorbing Brownian motion on $D$ and consequently almost excessive with respect to $X_D$ for each $\xi \in \partial D$. Therefore Lemma 7.1 means that the Martin kernel $K(x,\xi)$ is a Poisson kernel with respect to $\mu$ in the sense of Section 2. Hence, by defining the Feller kernel as (2.13) in terms of the present Martin kernel, we have the expression of the Feller measure

$$(7.6) \qquad U(d\xi, d\eta) = U(\xi,\eta)\mu(d\xi)\mu(d\eta).$$

We also see by Lemma 4.2 that $\mu$ is an admissible measure for $\partial D$ in the sense of Section 4.

On the other hand, we can see from (7.2) and (6.6) that the energy measure $\mu_{\langle u \rangle}$ of $u \in \mathcal{F}_e$ admits the expression

$$(7.7) \qquad \mu_{\langle u \rangle}(dx) = |\nabla u|^2(x)\,dx,$$

which does not charge the boundary $\partial D$. Hence all the conditions of Theorem 6.2 are satisfied for $F = \partial D$.

THEOREM 7.1. (i) *The measure $\mu$ on $\partial D$ defined by* (7.1) *is admissible with respect to the form* (7.2) *in the sense of Section* 4.

(ii) *For any $\mathcal{E}$-quasicontinuous $u \in L^2(\partial D; \mu) \cap \mathcal{F}_e$,*

$$(7.8) \qquad \mathcal{E}(Hu, Hu) = \tfrac{1}{2} \int_{\partial D \times \partial D} (u(\xi) - u(\eta))^2 U(\xi,\eta)\mu(d\xi)\mu(d\eta),$$

*where $Hu(x) = E^x(u(X_T); T < \infty), x \in \overline{D}$, and $U(\xi,\eta)$ is the Feller kernel defined in terms of the Martin kernel $K$.*

(iii) *Let $Y$ be the time-changed process of $X$ by means of PCAF with Revuz measure $\mu$. Then $Y$ is recurrent and of pure jump. In addition, $Y$ admits as its Lévy system*

$$(7.9) \qquad (U(\xi,\eta)\mu(d\eta), t),$$

*where $t$ denotes the nonrandom PCAF $\phi(t) = t$ of $Y$.*

(iv) *Let $(\check{\mathcal{F}}, \check{\mathcal{E}})$ be the Dirichlet space on $L^2(\partial D, \mu)$ of the time-changed process $Y$. Then*

(7.10)
$$\check{\mathcal{F}} = \left\{ f \in L^2(\partial D; \mu) : \int_{\partial D \times \partial D} (f(\xi) - f(\eta))^2 \times U(\xi, \eta) \mu(d\xi) \mu(d\eta) < \infty \right\},$$

(7.11) $\check{\mathcal{E}}(f, f) = \frac{1}{2} \int_{\partial D \times \partial D} (f(\xi) - f(\eta))^2 U(\xi, \eta) \mu(d\xi) \mu(d\eta), \qquad f \in \check{\mathcal{F}}.$

PROOF. Part (ii) follows from Theorem 6.2 and (7.6). Part (iii) follows from (ii) and Corollary 5.1. As for (iv), the inclusion $\subset$ in (7.10) and identity (7.11) are clear from (ii) and (6.1). Suppose that a function $f$ belongs to the space that appears in the right-hand side of (7.10). By virtue of [4], Theorem 3.1, we then have the expression of the function $w(x) = Hf^2(x) - (Hf(x))^2$,

$$w(x) = R^0_{0-} |\nabla(Hf)|^2(x), \qquad x \in D,$$

where $R^0_\alpha$ denotes the resolvent operator of the absorbing Brownian motion on $D$. Hence by setting $H^\alpha 1(x) = \int_{\partial D} K_\alpha(x, \xi) \mu(d\xi)$ by the kernel defined in (2.13), we easily see that

$$\int_D |\nabla(Hf)(x)|^2 \, dx = \lim_{\alpha \to \infty} \alpha (H^\alpha 1, w)_D.$$

From identity (6.14) we see that the right-hand side of the above equality is dominated by

$$\int_{\partial D \times \partial D} (f(\xi) - f(\eta))^2 \, U(\xi, \eta) \mu(d\xi) \mu(d\eta) < \infty,$$

proving that $Hf \in H^1_e(D)$ (see [11], Example 1.6.1) and consequently $f \in \check{\mathcal{F}}$. □

Equations (7.8) and (7.10) recover the Douglas integral description of the space of harmonic functions with finite Dirichlet integrals in [4] (but with the Feller kernel instead of the Naim kernel) for the present specific Martin space (cf. [8]).

**8. Reduction to Hunt processes.** This section is devoted to the proof of the following general reduction theorem especially applicable to the time changed process $Y$ in Section 4.

THEOREM 8.1. *Let $(E, m, \mathcal{F}, \mathcal{E})$ be a regular Dirichlet space and let $X = (X_t, P^x)$ be a right process over a subset $E_1 \subset E$ with $\operatorname{Cap}(E \setminus E_1) = 0$.*

We assume that $X$ is properly associated with $\mathcal{E}$ in the sense that $p_t u$ is an $\mathcal{E}$-quasicontinuous version of $T_t u$ for any $u \in L^2(X;m)$, where $p_t$ (resp. $T_t$) is the transition function of $X$ [resp. the $L^2$ semigroup associated with $(\mathcal{E}, \mathcal{F})$]. We further assume that the left limit $X_{t-}$ exists in $E_\Delta$ for every $t > 0$. Then there exists a Borel set $E_2 \subset E_1$ such that $\mathrm{Cap}(E \setminus E_2) = 0$, $E_2$ is $X$-invariant and the restriction $X|_{E_2}$ of $X$ to $E_2$ is a Hunt process properly associated with $\mathcal{E}$.

We prepare two lemmas.

LEMMA 8.1. (i) *For an open set $A \subset E$ of finite capacity, the function*
$$p_A^1(x) = E^x(\exp^{-\sigma_A}), \qquad x \in E_1,$$
*is an $\mathcal{E}$-quasicontinuous version of the $1$-equilibrium potential $e_A \in \mathcal{F}$ of $A$. Here $\sigma_A$ denotes the hitting time of the process $X$ for the set $A$.*

(ii) *If $\{A_n\}$ is a decreasing sequence of open subsets of $E$ with $\lim_n \mathrm{Cap}(A_n) = 0$, then*
$$\lim_{n \to \infty} p_{A_n}^1(x) = 0 \qquad \text{for } \mathcal{E}\text{-q.e. } x \in E_1.$$

PROOF. (i) It is known that $p_A^1$ is a version of $e_A$ (cf. [11], Lemma 4.2.1). Since $p_t p_A^1$ is an $\mathcal{E}$-quasicontinuous version of $T_t e_A$, we get the result by letting $t \downarrow 0$.

(ii) Since $\mathcal{E}(e_{A_n}, e_{A_n}) \downarrow 0$ as $n \to \infty$, (ii) follows from (i). □

LEMMA 8.2. *For any set $N \subset E_1$ with $\mathrm{Cap}(N) = 0$, there exists a Borel set $E' \subset E_1 \setminus N$ such that $\mathrm{Cap}(E \setminus E') = 0$ and $E'$ is $X$-invariant:*
$$P^x(X_t \in E'_\Delta \text{ for all } t \geq 0, \ X_{t-} \in E'_\Delta \text{ for all } t > 0) = 1,$$
*for all $x \in E'$.*

PROOF. There is a decreasing sequence of open sets $A_n$ including the set $(E \setminus E_1) \cup N$ such that $\lim_{n \to \infty} \mathrm{Cap}(A_n) = 0$. Lemma 8.1 then implies that
$$P^x(X_t \text{ or } X_{t-} \in B_0 \text{ for some } t \geq 0) = 0 \qquad \forall\, x \in E_1 \setminus N_1,$$
where $B_0 = \bigcap_n A_n \ [\supset (E \setminus E_1) \cup N]$ and $N_1$ is some subset of $E_1$ with $\mathrm{Cap}(N_1) = 0$.

Next we find a decreasing sequence of open sets $A'_n \supset B_0 \cup N_1$ with $\lim_{n \to \infty} \mathrm{Cap}(A'_n) = 0$ and let $B_1 = \bigcap_n A'_n$. Repeating the same argument, we can find an increasing sequence of Borel subsets $\{B_n\}$ of zero $\mathcal{E}$-capacity containing $(E \setminus E_1) \cup N$ such that
$$P^x(X_t \text{ or } X_{t-} \in B_n \text{ for some } t \geq 0) = 0 \qquad \text{for all } x \in E \setminus B_{n+1}.$$

Put $B = \bigcup_n B_n$. Then $E' = E \setminus B$ satisfies the desired properties. $\square$

PROOF OF THEOREM 8.1. From Lemma 8.1, we can see as in the proof of [11], Lemma 4.2.2, that for any $\mathcal{E}$-quasicontinuous function $u$ on $E$,

$$P^x\left(\lim_{t'\uparrow t} u(X_{t'}) = u(X_{t-}) \ \forall\, t > 0\right) = 1, \qquad \mathcal{E}\text{-q.e. } x \in E_1.$$

Choose a countable subfamily $C_1$ of $\mathcal{F} \cap C_0(E)$ which is dense in $C_0(E)$ and denote by $Q^+$ the set of all positive rational numbers. Since the functions $p_s f$ for $s \in Q^+$, $f \in C_1$ are $\mathcal{E}$-quasicontinuous, we can find a set $N$ with $\mathrm{Cap}(N) = 0$ such that the above identity holds for each $u = p_t f$, $s \in Q^+$, $f \in C_1$ and for all $x \in E_1 \setminus N$. We then use Lemma 8.2 to get a Borel set $E_2 \subset E_1 \setminus N$ such that $E_2$ is $X$-invariant and $\mathrm{Cap}(E \setminus E_2) = 0$. Since $X|_{E_2}$ is a right process on $E_2$ and

$$P^x\left(\lim_{t'\uparrow t} p_s f(X_{t'}) = p_s f(X_{t-}) \ \forall\, t > 0\right) = 1$$

for all $x \in E_2$ and for any $s \in Q^+$, $f \in C_1$, we can also prove that $X|_{E_2}$ is quasi-left continuous on $[0, \infty)$ in exactly the same manner as in the proof in [11], Lemma 7.2.5, completing the proof that $X|_{E_2}$ is a Hunt process on $E_2$. $\square$

M. FUKUSHIMA  
DEPARTMENT OF MATHEMATICS  
KANSAI UNIVERSITY  
OSAKA  
JAPAN  
E-MAIL: fuku@ipcku.kansai-u.ac.jp

P. HE  
DEPARTMENT OF APPLIED MATHEMATICS  
SHANGHAI UNIVERSITY OF FINANCE  
AND ECONOMICS  
SHANGHAI  
CHINA

J. YING  
INSTITUTE OF MATHEMATICS  
FUDAN UNIVERSITY  
SHANGHAI  
CHINA  
E-MAIL: jgying@fudan.edu.cn